\documentclass[12pt]{iopart}

\usepackage{iopams}
\usepackage{amssymb}
\usepackage{amsthm}
\usepackage{amscd}
\usepackage{epsfig}
\usepackage{graphicx}
\usepackage{subfigure}
\eqnobysec
\newtheorem{theorem}{Theorem} [section]        % 1st argument is your name for it
\newtheorem{proposition}[theorem]{Proposition} % 2nd argument is what is printed
\newtheorem{corollary}[theorem]{Corollary}
\newtheorem{lemma}[theorem]{Lemma}
\newtheorem*{remark}{Remarks}\newtheorem*{remarke}{Remark}

\normalsize
\def\N{{\rm I\kern-2pt N}}
\def\I{{\rm 1\kern-3.4pt l}}
\def\R{{\rm I\kern-2pt R}}
\def\P{{\rm I\kern-2pt P}}
\def\G{{\rm I\kern-2pt G}}
\def\E{{\rm I\kern-2pt E}}
\def\D{{\rm I\kern-2pt D}}
\def\C{{\rm I\kern-6pt C}}
\def\Z{{\rm Z\kern-4pt Z}}
\def\1{{\rm I\kern-6pt 1}}
\def\2{{\rm I\kern-4.6pt 1}}

\newcommand{\Dim}{{\rm Dim\,}}

\begin{document}

\title[{Phase transitions for self-similar measures}]% end with percent
 {Phase transitions for the   multifractal analysis of self-similar measures} % This is the full title of the paper
% Avoid equations in title, but if you insist: $E=\lowercase{mc}^2$
% Do not use the \thanks{} command; use \extraline{} instead (see above).

\author{B Testud}

\address{MAPMO UMR 6628, F\'ed\'eration Denis Poisson,  Universit\'e d'Orl\'eans, BP 6759, 45067 Orl\'eans cedex 2, France}
\ead{benoit.testud@univ-orleans.fr}
\begin{abstract}
We are interested  in   the  multifractal analysis  of a class of self-similar measures with overlaps. This class,  for which we obtain explicit formulae for  the $L^q$-spectrum $\tau(q)$ as well as the singularity spectrum $f(\alpha)$, is sufficiently large to   point out  new phenomena in the multifractal structure of self-similar  measures. We show that, unlike the  classical quasi-Bernoulli case,  the $L^q$-spectrum $\tau(q)$ of the measures studied  can have an   arbitrarily large number   of non-differentiability points (phase transitions).  These singularities occur only    for the  negative values of  $q$ and  yield to measures that do not satisfy  the usual  multifractal formalism. The  weak quasi-Bernoulli property is the key point of most of   the arguments. 
\end{abstract}
\submitto{\NL}
%Uncomment for PACS numbers title message
%\pacs{00.00, 20.00, 42.10}
% Keywords required only for MST, PB, PMB, PM, JOA, JOB? 
%\vspace{2pc}
%\noindent{\it Keywords}: Article preparation, IOP journals
% Uncomment for Submitted to journal title message
%\submitto{\JPA}
% Comment out if separate title page not required
%\maketitle

\section{Introduction}

Let us begin with some notation. For an integer $\ell\ge 2$,  we denote by  $\mathcal{F}=\cup_n \mathcal{F}_n$ where $\mathcal{F}_n$ is the set of the 
$\ell$-adic intervals of the $nth$ generation   included in the interval $\lbrack 0,1)$.  In other terms,
$\mathcal{F}_n =\left\{ I=\lbrack k/\ell^n,(k+1)/\ell^n)\, ,\, 0\le
k<\ell^n\right\}.$ For every $x\in \lbrack 0,1)$, $I_n(x)$ stands for the unique interval among $\mathcal{F}_n$ containing~$x$.

Let $m$ be  a  probability measure on the interval  $\lbrack 0,1)$. For  $x\in \lbrack 0,1)$, we  define the {\em local dimension} (also called H\"older exponent) of $m$ at $x$ by \[\alpha(x)= \lim_{n\rightarrow +\infty} -{{\log m (I_n(x))}\over {n \log \ell}},\] provided this limit exists.  The  aim of  multifractal analysis is  to find the Hausdorff dimension, $\dim(E_\alpha)$,  of  the  level  set  $ E_\alpha=\left\{x : \alpha(x)=\alpha\right\}$ for  $\alpha>0$.   The function  $f(\alpha)=\dim(E_\alpha)$ is called the {\em singularity spectrum} (or multifractal spectrum) of $m$ and we say that $m$ is a {\em multifractal measure} when $f(\alpha)>0$ for several $\alpha{}'s$.

The concepts underlying the  multifractal decomposition of a measure  go back to an early  paper of    Mandelbrot \cite{mandel}.   In the 80's multifractal measures were used    by physicists to study various  models arising from natural phenomena.   In    fully developped turbulence they  were  used by Frisch and Parisi \cite{fp} to investigate the intermittent behaviour in the regions of high vorticity. In dynamical system theory they were used by Benzi et al.~\cite{ben}  to measure how often a given region of the attractor is visited.  In diffusion-limited aggregation (DLA) they were  used by Meakin et al.~\cite{mea}  to describe the probability of a random walk landing to the neighborhood of a given site on the aggregate. 
   
In order to determine the function $f(\alpha)$,       Hentschel and  Procaccia  \cite{hp}  used  ideas based on  Renyi entropies \cite{re} to  introduce the generalized dimensions  $D_q$ defined by  \[D_q= \lim_{n\rightarrow +\infty} {1\over {q-1} } {{\log\left(\sum_{I\in\mathcal{F}_n} m(I)^q\right)}\over {n \log \ell}},\] (see also \cite{grp,gr}). From a physical and heuristical point of view, Halsey et al. \cite{has} showed that the singularity spectrum $f(\alpha)$ and the generalized dimensions $D_q$ can be derived from each other. The Legendre transform  turned out to be a useful tool linking  $f(\alpha)$ and $D_q$.  More precisely,  it was suggested that   
\begin{eqnarray}\label{fm} f(\alpha)=\dim(E_{\alpha})=\tau^*(\alpha)=\inf(\alpha q+\tau(q),\ q\in \R),\end{eqnarray}  where \[\tau (q) = \limsup_{n\rightarrow +\infty}\,\tau _n (q)\quad\mbox{ with }\quad \tau _n (q)={1\over {n\log \ell}}\log\left(\sum_{I\in\mathcal{F}_n} m(I)^q\right).\]   (The sum  runs over the $\ell$-adic intervals $I$  such that $m(I)\not=0$.)   The  function $\tau(q)$ is  called the  $L^q$-spectrum of $m$   and if the limit exists $\tau(q)=(q-1)D_q$. Note that there may be  problems   of stability or  invariance   in the definition  of $\tau(q)$ for negative $q$   and   Riedi  \cite{ri} propose  an improvement of this definition.     In what follows,  these difficulties   will be  avoided by restricting the sums over   convenient $\ell$-adic intervals defining  the  measure. Of course, this way  is not an option in many applications where the structure of the measure is not known in advance. For more information on  the $L^q$-spectrum and the singularity spectrum we refer  the reader to    \cite{bn,bnbh,bmp,fan,feng,h,ngg,o1,p,t1,t3,t2,yyl}.

Relation \eref{fm} is called the multifractal formalism and in many aspects it is analogous to the well-known thermodynamic formalism developed by Bowen \cite{bo} and Ruelle \cite{ru}. 

 For number of measures, relation \eref{fm}    can  be verified rigorously. In particular, under some separation conditions, 
   self-similar measures satisfy the multifractal formalism (e.g. \cite{cm,em,fal,lng,o}).  Despite all the investigations mentioned, the exact range of the   validity of  the multifractal formalism is still not known. Furthermore, it is easy to construct measures that do not  satisfy     \eref{fm}  (e.g \cite{ri}). It  is   thus interesting to find conditions  ensuring  the validity of  \eref{fm}. The main difficulty is often to get  a lower bound of $\dim(E_\alpha)$.  Usually, such a minoration relies on the
existence of an auxiliary measure $m_q$, the so-called    Gibbs measure,   supported on  the level set $E_\alpha$. Recall  that  $m_q$ is a {\em Gibbs measure} at  state $q$ for the measure $m$  if 
\begin{eqnarray*}\forall n,\,\,\forall I\in \mathcal{F}_n, \quad
  {1\over C}  m(I)^q \ell^{-n\tau(q)}\le  m_q(I)\le C m(I)^q
  \ell^{-n\tau(q)},\end{eqnarray*} where the constant  $C>0$ is
independent of $n$ and $I$. If  $\tau$ is differentiable at $q$,   the measure  $m_q$ --if it exists--  will be 
supported by   $E_{-\tau'(q)}$. In this case,     Brown, Michon
and Peyri\`ere    established  \cite{bmp,p}  that
\begin{eqnarray*}
\dim(E_{-\tau'(q)})
=\tau^*(-\tau'(q))=-q\tau'(q)+\tau(q).\end{eqnarray*} 
In general, to prove the existence of Gibbs measures  we need some   homogeneity hypotheses on the measure.  This is, for instance, the case of     quasi-Bernoulli measures~:    a probability measure $m$ is said to be  {\em quasi-Bernoulli} if there exists a constant $C>0$ such that 
\begin{eqnarray}\label{qb}\fl \forall (n,p)\in\N^2,\,\forall I\in \mathcal{F}_n,\,\forall J\in \mathcal{F}_p,\quad {1\over C} m(I)m(J)\le m(IJ)\le C m(I)m(J),\end{eqnarray}
where  $IJ=I\cap\sigma^{-n}(J)$ and $\sigma(x)=\ell x\,(\mbox{mod}\,1)$
is the  shift map on the interval $[0,1)$. In this situation,  Brown, Michon
and 
Peyri\`ere  \cite{bmp,m,p}  proved     the existence of a Gibbs measure at
every state $q$. A few years later,    Heurteaux   \cite{h}  showed  that $\tau$ is differentiable on $\R$. Therefore, for   quasi-Bernoulli measures,  we have \[\forall \alpha  \in \left( -\tau'(+\infty), -\tau'(-\infty) \right) ,\quad
\dim(E_\alpha)=\tau^*(\alpha).\]

Recently, in \cite{t1,t2}  we introduced  a more general condition that we call   the weak quasi-Bernoulli
property. More precisely, we say that   a measure   $m$ satisfies   the {\em weak  quasi-Bernoulli} property  if there exists  a constant $C>0$
and some   integers $r_1, r_2, p_1, p_2, s_1, s_2$ such that \begin{eqnarray}\label{wqb}\fl
\cases{\exists C>0,\,  \forall n,\forall p ,\, \forall I\in \mathcal{F}_n, \forall J\in \mathcal{F}_p,&\\C^{-1} m(I) \displaystyle\sum_{k=r_1}^{r_2} m(\sigma^{-k}(J))\le   \displaystyle\sum_{k=p_1}^{p_2} m(I\cap\sigma^{-(n+k)}(J))\le C m(I)  \displaystyle\sum_{k=s_1}^{s_2} m(\sigma^{-k}(J)).&\\}
\end{eqnarray}  At first sight, this new condition may seem artificial but is in fact 
    natural. Indeed, in \cite{t1,t2} we showed that     many  self-similar measures
  with overlaps    are  not
    quasi-Bernoulli but are weak quasi-Bernoulli and  may be used to estimate
    the dimension of  self-affine graphs.

  Furthermore, under this condition,  we proved in \cite{t2} the existence of  Gibbs measures  at
 every positive state $q$ and the differentiability of   $\tau$  on $\R^+$. For  weak quasi-Bernoulli measures, we 
 deduced that  
 \begin{eqnarray*}\forall \alpha  \in \left( -\tau'(+\infty),
     -\tau'(0) \right),\quad
   \dim(E_\alpha)=\tau^*(\alpha).\end{eqnarray*}

Now, it is natural to ask whether or not  these  results 
 still hold  for  negative
 $q$ when the measure  only  satisfies the weak quasi-Bernoulli
property. In  particular, in this setting, we would like to know if \begin{enumerate}\item  the $L^q$-spectrum
  $\tau(q)$ is  differentiable on  $(-\infty,0)$, \item there exists   Gibbs
  measures for negative  $q$,      \item we have
  $\dim(E_\alpha)=\tau^*(\alpha)$ for  $\alpha>-\tau'(0)$.\end{enumerate}

Note that  the tools   used  in this context for  $q\ge 0$ cannot be applied for $q<0$. In particular, to prove the existence of   Gibbs measures for  $q\ge 0$   we use some 
multiplicative properties 
of the sequence $\ell^{n\tau_n(q)}$ which are  no longer verified for   $q<0$. In what follows, we   show that the multifractal formalism may  break down for weak quasi-Bernoulli measures.  Therefore, for these
measures,  the  answer to the above  questions could  be  no. 

%This class    may appear rather restrictive but is in fact sufficiently large to observe new and interesting phenomena in  the multifractal structure of self-similar measures. 

Let us  precise these examples.  For an integer $\ell\ge 2$,  we
consider the $2 \ell$  similitudes $S_i: [0,1]\mapsto
[0,1]$  defined by \[\forall\, 0\le i \le \ell-1, \quad S_i(x)={1\over
  {\ell}}x +{{i}\over {\ell}}\quad \mbox{and} \quad S_{i+\ell}(x)=-{1\over
  {\ell}}x +{{i+1}\over {\ell}}.\] For a given probability weight
$\{p_i\}_{i=0}^{2\ell-1}$,  it is well known (e.g \cite{ fa, hu}) that there exists
a unique probability   measure $\mu$ on $[0,1]$   verifying
\begin{eqnarray} \label{mu} \mu=\sum_{i=0}^{2\ell-1} p_i\,\mu\circ
  S_i^{-1}.\end{eqnarray}This measure is often called the self-similar measure
generated by $\{S_i\}_{i=0}^{2\ell-1}$. In this paper we establish  that $\mu$
satisfies the weak quasi-Bernoulli property.   Moreover, we show that          there exists a {\em Frostman measure}
$\mu_q$ at every negative state $q$, i.e. a  measure  $\mu_q$ such that
\begin{eqnarray}\label{fr}\forall n,\,\,\forall I\in \mathcal{F}_n,  \quad
  \mu_q(I)\le C \mu(I)^q \ell^{-n\tau(q)},\end{eqnarray} where the constant $C>0$  is independent of $n$ and $I$. Thus, for  $\alpha=-\tau'(q)$,   we have 
   \begin{eqnarray*}\forall x\in E_\alpha,  \quad \mu_q(I_n(x))\le (\ell^{-n})^{\tau^*(\alpha)},\end{eqnarray*} if   $n$ is   large enough. The mass distribution principle or  Frostman Lemma  (e.g. \cite{fa})  implies that $\dim(E_\alpha)\ge \tau^*(\alpha)$. Thus, the   values of $\alpha$ for which the
 multifractal formalism  may fail  lie   in   intervals   
 $(-\tau'_+(q),-\tau'_-(q))$ where $q$ is a point of non-differentiability of $\tau$  ($\tau'_-(q)$ and $\tau'_+(q)$ stand for the left and the right derivative  respectively).  Such a  point $q$  will be called   a  {\em phase transition}.

 We assume  that the weights $p_i$ associated to the measure $\mu$ verifying \eref{mu}  are positive  for every   $0\le i\le \ell-1$. We set $B=\left\{0\le i\le\ell-1,\,\,p_{i+\ell}=0\right\}$ and $\tilde\tau(q)= \log_\ell\left(\sum_{i\in B}
p_i{}^q\right)$.   In this case, the $L^q$-spectrum $\tau_\mu(q)$ of  $\mu$ is given by $\tau_{\mu}(q)=\max\left(\tau_{\nu}(q),\tilde\tau(q)\right)$ where  $\nu=(\mu+\mu\circ T)/2$  and $T(x)=1-x$.  In order to get phase transitions for the function $\tau_\mu$, it is thus enough to find conditions   on  the  $p_i\,'s$ ensuring that  the equation
$\tau_\nu(q)=\tilde\tau(q)$ has  isolated solutions.

Let $K$ be the compact  set defined by
$K=\cup_{i\in B} S_i(K)$. The attractor $K$ plays an important   role to determine   the local dimensions  of $\mu$. Indeed,  we can link the level sets of   $\mu$ and $\nu$ in the following way~:  $E_{\alpha}(\mu)=\left(E_{\alpha}(\nu)\cap
([0,1]\setminus K) \right)\cup\left(E_{\alpha}(\mu)\cap
K\right).$  If $\nu$ satisfies
the quasi-Bernoulli property, we get $\dim\left(E_\alpha(\mu)\right)=\max\left(\tau_\nu{}^*(\alpha),
\tilde\tau^*(\alpha)\right).$  Using  the expression of $\tau_\mu$, we   deduce that each phase
transition   corresponds   to an  interval  in which   the multifractal formalism does not hold. More precisely,  we have the following.      \begin{enumerate} \item If
  $\tau_{\mu}{}'(q)$ exists and if   
  $\alpha=-\tau_{\mu}{}'(q)$, then  
  $\dim(E_{\alpha}(\mu))=\tau_\mu{}^*(\alpha).$\item If
  $\tau_{\mu}{}'(q)$ does not exist and if
  $-(\tau_{\mu})_+'(q)<\alpha<-(\tau_{\mu})_-'(q),$ then   $\dim(E_{\alpha}(\mu))<\tau_\mu{}^*(\alpha).$\end{enumerate} The class of  measures studied    may appear rather restrictive but is in fact sufficiently large to point out  new and interesting phenomena.  In particular, we
can observe the following facts. 
    \begin{itemize}\item  The existence of an
  isolated point in the set of the local dimensions $D_\mu$  defined by   $D_\mu=\{\alpha,E_\alpha(\mu)\not=\emptyset\}$. This situation has  already been obtained for the Erd\"os measure 
   and for the   $3$-time convolution of the Cantor measure 
  (e.g \cite{fl1,hl}).\item The existence of non-concave multifractal   spectra eventually
  supported by a  union of mutually  disjoint  intervals. To the best of our 
  knowledge, it is the first time that such multifractal structures are
  obtained for  self-similar measures. \item The existence of   an arbitrarily large number of phase  transitions for the $L^q$-spectrum $\tau(q)$.\end{itemize}

The paper is organized as follows. In section~2 we prove that the measure $\mu$, given by \eref{mu}, satisfies the weak quasi-Bernoulli property. In  section~3 we establish   the existence of Gibbs measures at every negative state $q$ for the measure $\mu$.  In  section~4 we determine the $L^q$-spectrum $\tau_\mu(q)$. In section~5 we are interested in the singularity spectrum of $\mu$.    The paper ends with a range of examples. 

\section{The weak quasi-Bernoulli property}

Let us introduce some notation. In what follows, except contrary mention,  $\ell\ge 2$ is  an integer, $I$   a $\ell$-adic
interval of $nth$ generation and $J$    a $\ell$-adic interval.  For every  $(\epsilon_1 ,...,\epsilon_n) \in \{0,...,\ell-1\}^n$,
$I_{\epsilon_1\cdots \epsilon_n}$ stands for the element of $\mathcal{F}_n$
defined by  \[I_{\epsilon_1\cdots
  \epsilon_n}=\left\lbrack\sum_{i=1}^n{{\epsilon_i}\over {\ell^i}},
  \sum_{i=1}^n{{\epsilon_i}\over {\ell^i}}+{1\over{\ell^n}}\right).\] If
$I=I_{\epsilon_1\cdots\epsilon_n}$ and $\epsilon\in\{0,\cdots,\ell-1\}$, we
shall write $\epsilon I$  instead of
$I_{\epsilon\epsilon_1\cdots\epsilon_n}$. If $f$ and $g$ are positive functions of the same parameter,  $f\approx g$ means  there exists a constant $C>0$ such that $C^{-1}g\le f \le C g$. Moreover, for any  matrices  $M$ and $N$,   we shall write $M>0$ (and we shall say that $M$ is positive) if all the digits of $M$ are positive and $M>N$ if $M-N>0$. The matrice relations $<$, $\ge$ and $\le$ are similarly  defined. Finally, for a  $2\times 2$ nonnegative matrix $M$, we define    $\Vert M\Vert_1$ and $\Vert M\Vert$ by 
 \[\Vert M\Vert_1=\left(\begin{array}{cc} 1 & 0\end{array}\right) M \left(\begin{array}{c} 1\\
   1\end{array}\right) \quad \mbox{and} \quad \Vert M\Vert={1\over 2}
 \left(\begin{array}{cc} 1 & 1\end{array}\right) M \left(\begin{array}{c} 1\\ 1\end{array}\right).\]

Let $\mu$ be the measure verifying \eref{mu}. For convenience,  we suppose that $\mu$ is supported on the interval $[0,1]$. That is
equivalent to the condition~: $p_i+p_{i+\ell}>0$, for every $i\in\{0,\cdots,\ell-1\}$. The relation \eref{mu} implies that  \begin{eqnarray*}\forall  \epsilon \in\{0,\cdots,\ell-1\},\quad     \mu(\epsilon I)=p_{\epsilon}\,\mu(I)+p_{\epsilon+\ell}\,\mu(I^*),\end{eqnarray*} where $I^*=T(I)$ and $T(x)=1-x$. Since $(I^*)^*=I$,  we have
\begin{eqnarray}\label{mat}\fl\left(\begin{array}{c}\mu(\epsilon I) \\  \mu\circ T(\epsilon I)\end{array}\right) = M_{\epsilon}\left(\begin{array}{c}  \mu(I) \\  \mu\circ T(I)\end{array}\right) \quad\mbox{where} \quad   M_{\epsilon}=\left(\begin{array}{cc} p_{\epsilon} & p_{\epsilon+\ell} \\p_{2\ell-1-\epsilon}&p_{\ell-1-\epsilon}\end{array}\right).\end{eqnarray}     By iterating this relation,   we get \begin{eqnarray*}\forall I= I_{\epsilon_1\cdots \epsilon_n}\in \mathcal{F}_n,\quad\left(\begin{array}{c}  \mu(I) \\\mu\circ T(I)\end{array}\right) = M_{\epsilon_1}\cdots M_{\epsilon_n}\left(\begin{array}{c}  1 \\1\end{array}\right), \end{eqnarray*}  and we   deduce that 
 \begin{eqnarray}\label{norme} \mu(I)=\Vert M_I\Vert_1\quad \mbox{and}\quad
   \nu(I)=\Vert M_I\Vert,  \quad\mbox{where} \quad   M_I= M_{\epsilon_1}\cdots M_{\epsilon_n}.
\end{eqnarray}  
These relations will  be used to prove that     $\mu$  satisfies  the  weak quasi-Bernoulli property~\eref{wqb}. More precisely, we have the following.
\begin{proposition}\label{prop1} Let  $\mu$ be the measure verifying  \eref{mu}. Then,\[
\cases{ \exists C>0,\,  \forall n,\forall p ,\, \forall I\in \mathcal{F}_n, \forall J\in \mathcal{F}_p,&\\ C^{-1} \mu(I)\mu(J)\le \mu(I\cap \sigma^{-(n+1)}(J))\le C \mu(I)\mu(\sigma^{-2}(J)),&\\}
\] where $\sigma :[0,1]\mapsto [0,1]$ is the shift map on the $\ell$-adic basis given by $\sigma(x)=\ell x (mod 1)$.  \end{proposition}
\begin{remark} \rm \begin{enumerate} \item[1.] In general,   the measure $\mu$ does not satisfy the
  quasi-Bernoulli property~\eref{qb}. To see that, take for example   $\ell=2$ and suppose that  $p_0>p_1$,  $p_0p_1p_2>0$ and  $p_3=0$.   Using \eref{mu}, we get for     $J=I_{1\cdots1}\in\mathcal{F}_n$, \[\mu(0J)=\mu(I_{01\cdots1})=p_0
  \mu(J)+p_2\mu(J^*)=p_0\mu(I_{1\cdots1})+p_2\mu(I_{0\cdots0}).\] From  \eref{mat} and \eref{norme},
     $\mu(J)=\Vert M_1{}^n\Vert_1=p_1{}^n$ and $\mu(J^*)=\Vert M_0{}^n\Vert_1 \approx
  p_0{}^n$. Therefore,  if  $I=I_0$,  we have   $\mu(IJ)\approx p_0{}^n$ and   $\mu(I)\mu(J)\approx p_1{}^n$, which proves   that   $\mu$   is not quasi-Bernoulli.

\item[2.] If  for every 
  $\epsilon\in\{0,\cdots,\ell-1\}$,  $p_\epsilon p_{\epsilon+\ell}=0$,  the Open Set
  Condition of Hutchinson~\cite{hu} is verified. In this case,  $\mu$ is quasi-Bernoulli and proposition \ref{prop1} easily follows.\end{enumerate}\end{remark}

\smallskip

\noindent{\bf Proof of  proposition~\ref{prop1}.}  According to the above remark, we can suppose that there exists  $\tilde\epsilon\in\{0,\cdots,\ell-1\}$ such that   $p_{\tilde\epsilon} p_{\tilde\epsilon+\ell}>0$. By \eref{mat},    $M_{\tilde\epsilon}+M_{\ell-1-\tilde\epsilon}>0$. Thus,   we can find  a constant $C>0$ such that
\[{1\over C}E\le \sum_{\epsilon=0}^{\ell-1} M_\epsilon \,\,\,\,\mbox{ and }\,\,\,\, I_2\le C E\sum_{\epsilon=0}^{\ell-1}M_\epsilon,\]
where \[E=\left(\begin{array}{cc} 1&0\\1&0\end{array}\right)\,\,\,\,\rm{ and }\,\,\,\,I_2=\left(\begin{array}{cc} 1&0\\0&1\end{array}\right).\] 
It follows from  \eref{norme} that  \[\fl\mu(I\cap \sigma^{-(n+1)}(J))=\left\Vert M_I\left(\sum_{\epsilon=0}^{\ell-1}M_\epsilon\right)M_J\right\Vert_1\ge{1\over C}\Vert M_IEM_J\Vert_1={1\over C}\mu(I)\mu(J).\]
On the other hand, we obtain in a similar way that 
\[\fl\mu(I\cap \sigma^{-n}(J))=\Vert M_IM_J\Vert_1\le C \left\Vert M_IE\left(\sum_{\epsilon=0}^{\ell-1}M_\epsilon\right)M_J\right\Vert_1=C\mu(I)\mu(\sigma^{-1}(J)).\]
A monotone class  argument implies  that  this relation still
holds if we  replace $J$ by any Borel set  $B$. Thus,  by taking
$B=\sigma^{-1}(J)$, we have  
\[\mu(I\cap \sigma^{-(n+1)}(J))\le C \mu(I)\mu(\sigma^{-2}(J)),\]
which completes the proof of proposition \ref{prop1}.$\square$

\section{Frostman measures}
In this section we establish the existence of Frostman measures \eref{fr} at every negative state $q$ for the measure $\mu$ defined by \eref{mu}. We begin with  a preliminary result that  gives  conditions ensuring the existence of Frostman measures.

For a  probability measure $m$ on the
  interval $[0,1]$, let us define   the series $Z(s)$   by  \[\forall s\in \R, \quad Z(s)=\sum_{n\ge 1}
u_n\ell^{-ns}, \quad \mbox{where}\quad u_n=\ell^{n\tau_n(q)}.\]  If  $m(I)>0$,    $Z_I(s)$ denotes   the series associated to the measure $m_I$ verifying  $m_I(J)=m(IJ)/ m(I)$.

\begin{proposition}\label{prop2} Let  $m$ be a probability measure on the
  interval $[0,1]$ and  $q\in \R$.  With the above notation, suppose that       there exists a constant $C>0$ such that  \begin{enumerate}\item  $\forall n,\, \forall p,\quad u_{n+p}\le C u_n u_p$, \item $\forall I\in\mathcal{F},\, \forall s\in \R,\quad  Z_I(s)\le CZ(s)$\end{enumerate}
Then, there exists a Frostman measure at state $q$ for the measure $m$.\end{proposition} 

\noindent{\bf Proof.} We adapt to our situation  the arguments used by   Michon and Peyri\`ere    \cite{m,p} in another context.  The submultiplicativity property 
of the sequence $v_n=C\ell^{n\tau_n(q)}$ implies that the sequence $v_n^{1/n}$ tends to its lower bound. As a
consequence, $\tau_n(q)$ converges and if we call  $\tau(q)$ its limit, we
have 
\[\forall n\in\N,\quad C\ell^{n\tau_n(q)}\ge\ell^{n\tau(q)}.\]
Therefore, the series $Z(s)$ converges for  $s>\tau(q)$ and diverges for
$s=\tau(q)$. Let us consider, for $s>\tau(q)$, the function $\phi_s$ defined by
\[\phi_s(x)=\sum_{n\ge 1} m(I_n(x))^q\left(\ell^{-n}\right)^{-1+s}.\]Since   $\int_0^1 \phi_s(x) dx=Z(s)$, we can   define  a  probability measure  $\nu_s$ on the
interval $[0,1]$~by  \[\forall I\in\mathcal{F},\quad \nu_s(I)={{\int_I\phi_s(x)dx}\over
  {Z(s)}}.\] For every $I\in \mathcal{F}_n$, we find that    \[ Z(s)
\nu_s(I)=\int_I \phi_s(x)\, dx= \ell^{-n} \sum_{1\le k\le n} m(I_k)^q \left(\ell^{-k}\right)^{-1+s} +m(I)^q\ell^{-ns} Z_I(s),\] where $I_k$ denotes the element
of  $\mathcal{F}_k$   containing  $I$.

Let  $m_q$    be  a weak$^*$-limit  of
$\nu_s$  as $s$ goes to $\tau(q)$. The divergence of the series   $Z(s)$ for
$s=\tau(q)$ and the inequality    $Z_I(s)\le C Z(s)$ imply that 
$m_q(I)\le C m(I)^q\ell^{-n\tau(q)}$,  which completes the proof of 
proposition \ref{prop2}.$\square$ 

We   easily deduce the following result.
\begin{corollary} Let   $m$ be a probability measure on the interval $[0,1]$
  and $q\in \R$. Suppose that  there exists a constant  $C>0$ such that   \begin{eqnarray}\label{cp}\forall I,\,\forall J,\quad m(IJ)^q\le C m(I)^q m(J)^q.\end{eqnarray}
 Then, there exists   a Frostman measure at state $q$  for the measure $m$.

In particular, the condition \eref{cp} is satisfied if  $m(IJ)\le C m(I)m(J)$ and  $q>0$ or if $m(IJ)\ge C m(I)m(J)$ and $q<0$. \end{corollary}

 We will use the following lemma to prove  that the measure  $\mu$
 satisfies  the hypotheses of  proposition \ref{prop2}. 
\begin{lemma} \label{lemme1} Let $\mu$ be the measure defined by \eref{mu}. For every $I\in \mathcal{F}$,
   one of the following is satisfied. 
   \begin{enumerate}
 \item[(i)] $\forall J\in \mathcal{F},
 \quad  \mu(I)\mu(J)\le 2 \,\mu(IJ),$
  
 \item[(ii)] $\forall J\in \mathcal{F},
 \quad  \mu(I)\mu\circ T(J)\le 2\,\mu(IJ),$ where
$T(x)=1-x$.\end{enumerate}   \end{lemma} 
\noindent{\bf Proof.}  For 
 $I\in \mathcal{F}$ and $\epsilon\in\{0,\cdots,\ell-1\}$, we have   $S_i^{-1}(\epsilon I)=I$ or   $S_i^{-1}(\epsilon I)=T(I)$ or $S_i^{-1}(\epsilon
I)=\emptyset$.  Thus, by iterating \eref{mu}, we can 
find two non-negative  real numbers   $A(I)$ and $B(I)$ (depending only  on $I$) such that 
\begin{eqnarray}\label{imp}\forall J\in \mathcal{F},\quad \mu(IJ)=A(I)\mu(J)+B(I)\mu\circ T(J).\end{eqnarray}
  We then
obtain that either

\medskip
\hskip3.5cm $\forall J\in \mathcal{F}, \quad \mu(IJ)\ge C(I)\mu(J),$\\
or

\hskip3.5cm $\forall J\in \mathcal{F},
 \quad \mu(IJ)\ge  C(I)\mu\circ T(J),$

\medskip
\noindent
where   $C(I)=\max(A(I),B(I))>0$. 

By  taking $J=[0,1]$ in \eref{imp}, we  get   
  $\mu(I)\le 2C(I)$.   Hence, either

\medskip
\hskip3.5cm $\forall J\in \mathcal{F}, \quad 2 \mu(IJ)\ge \mu(I)\mu(J),$\\
or

\hskip3.5cm $\forall J\in \mathcal{F},
 \quad 2 \mu(IJ)\ge  \mu(I)\mu\circ T(J),$

\medskip
\noindent which completes the proof of   lemma \ref{lemme1}.$\square$ 

\begin{theorem} \label{fro}  There exists   a  Frostman measure  at every state $q<0$  for the measure $\mu$ verifying \eref{mu}.\end{theorem}

 \noindent {\bf Proof.}  Let  $I\in\mathcal{F}$. If for every    $J\in \mathcal{F}$,    $ \mu(I)\mu(J)\le 2 \mu(IJ)$ (respectively, $\mu(I)\mu\circ T(J)\le 2 \mu(IJ)$), 
we set    $\tilde{\mu_I}=\mu$ (respectively,
  $\tilde{\mu_I}=\mu\circ T$).   By lemma \ref{lemme1}, we have  
   \begin{eqnarray*} \forall J\in \mathcal{F},
 \quad  \mu(I)\tilde{\mu_I}(J)\le 2 \mu(IJ).\end{eqnarray*} Since $q$ is negative, we obtain that 
   \[\fl u_{n+p}=\sum_{I\in\mathcal{F}_n}
\sum_{J\in\mathcal{F}_p}\mu(IJ)^q\le \left({1\over 2}\right)^q   \sum_{I\in\mathcal{F}_n}
\sum_{J\in\mathcal{F}_p} \mu(I)^q\tilde{\mu_I}(J)^q=\left({1\over 2}\right)^q  u_n u_p\]
and \[\fl Z_I(s)=\sum_{n\ge 1}
    \sum_{J\in\mathcal{F}_n}\left({{\mu(IJ)}\over{\mu(I)}}\right)^q \ell^{-ns}
    \le \left({1\over 2}\right)^q \sum_{n\ge 1}\sum_{J\in\mathcal{F}_n}\tilde{\mu_I}(J)^q
    \ell^{-ns}= \left({1\over 2}\right)^q  Z(s).\]  Thus,   $\mu$ satisfies the hypotheses of proposition \ref{prop2} and     theorem \ref{fro} follows.$\square$

\section{The function $\tau_{\mu}$}
In this section we determine  the $L^q$-spectrum
$\tau_{\mu}(q)$ of the measure $\mu$ verifying \eref{mu}.   Let us start with the following easy lemma.
\begin{lemma}\label{tech}  Let  $(u_n)_{n\in\N}$ and  $(v_n)_{n\in\N}$ be two
  sequences of real positive  numbers  such that 
\[\lim_{n\rightarrow \infty}  u_n^{1/ n}=u \, \mbox{ et } \,
\lim_{n\rightarrow \infty}  v_n^{1/ n}=v.\]
Let   \[w_n=\sum_{k=0}^n u_k v_{n-k}.\] Then,  the sequence
$(w_n^{1/n})_{n\in\N}$   converges to  $w=\max\{u,v\}$.
\end{lemma} The proof is elementary and therefore omitted. 
\begin{theorem}\label{tau} Let  $\mu$  be the  probability
  measure verifying   \eref{mu}.     Suppose that for  every  $0\le i\le \ell-1$,
  $p_i>0$, and set    $B=\left\{0\le i\le \ell-1, p_{i+\ell}=0\right\}$.
 Then, by denoting  $\tilde\tau(q)=\log_\ell\left(\sum_{i\in B}p_i{}^q\right)$, we have 
\begin{eqnarray}\label{taut} \forall q \in \R,\quad
  \tau_\mu(q)=\max(\tau_{\nu}(q),\tilde\tau(q)),\end{eqnarray} where $\nu=(\mu+\mu\circ T)/2$ and $T(x)=1-x$. By convention, $\tilde\tau(q)=-\infty$  if $B=\emptyset$.
\end{theorem}

\begin{remark} \rm \begin{enumerate} \item[1.] We obtain a similar result  replacing the hypothesis ``$p_i>0$, for  every  $0\le i\le \ell-1$''  by ``$p_i+p_{i+\ell}>0$, for  every  $\,  0\le i\le \ell-1$''.
 
 \item[2.] For    positive  $q$, the two $L^q$-spectra $\tau_{\mu}(q)$ and
  $\tau_{\nu}(q)$  are obviously the
  same.  If $q<0$, the  situation is inverted~: ``small become big''. In fact, sets with negligeable mass determine the function $\tau_\mu$ ; therefore it suffices to    consider  the sum over  indices $i\in B$ in the expression of $\tilde\tau$.    

\item[3.] If $B=\{0,\cdots,\ell-1\}$, $\mu$ is a multinomial measure (also called
Bernoulli product). The calculation of $\tau_\mu$ is then straightforward and we have    $\tau_\mu=\tilde\tau$ (e.g  \cite{fal}).

\item[4.]  To obtain phase transitions for the function $\tau_\mu$, it is thus enough to find conditions   on  the   $p_i{}'s$ ensuring that  the equation
$\tau_\nu(q)=\tilde\tau(q)$ has many isolated  solutions.\end{enumerate}\end{remark}

\noindent{\bf Proof of  theorem~\ref{tau}.} We fix  $q\in\R$, $n\in\N^*$  and we write
$C$ for each    constant which  depends  on
$q$ but  not on $n$.
Using  \eref{mu},  we get 
\begin{eqnarray*} w_n&=&\sum_{I\in\mathcal{F}_n} \mu(I)^q=\sum_{i=0}^{\ell-1} \sum_{I\in\mathcal{F}_{n-1}} \mu(iI)^q\\
&=&\sum_{i\not\in B} \sum_{I\in\mathcal{F}_{n-1}}(p_i
\mu(I)+p_{i+\ell}\mu(I^*))^q+\sum_{i\in B}
\sum_{I\in\mathcal{F}_{n-1}}p_i{}^q \mu(I)^q\\ &\le&
C\sum_{I\in\mathcal{F}_{n-1}} \nu(I)^q+ \sum_{i\in B} p_i{}^q
w_{n-1}:=Cv_{n-1}+\sum_{i\in B} p_i{}^q w_{n-1}.\end{eqnarray*}
By induction,     we obtain
\[w_n\le C\sum_{k=0}^{n-1} \left(\sum_{i\in B} p_i{}^q\right)^k
  v_{n-(k+1)}:=C\sum_{k=0}^{n-1} u_k
  v_{n-(k+1)}.\]  We can find a minoration of the same type in a similar way. 

Furthermore,   the definition of $\Vert.\Vert$ in \eref{norme} 
implies that  $\nu(IJ)\le
2\nu(I)\nu(J)$.  Thus,  the sequence $v_n^{1/ n}$
converges and  theorem \ref{tau} easily follows    from  lemma \ref{tech}.$\square$

\section{The level sets   $E_\alpha$}

 In this section we link  the level sets $E_\alpha(\mu)$ and
$E_\alpha(\nu)$ associated  to the   measures $\mu$ and~$\nu$. 
\begin{proposition}\label{ealpha} The hypotheses are the same as    in theorem \ref{tau}. Let $K$ be the compact set defined by
  $K=\bigcup_{i\in B} S_i(K)$ with the convention that  $K=\emptyset$ if  $B=\emptyset$. Then,\[E_{\alpha}(\mu)=\left(E_{\alpha}(\nu)\cap
([0,1]\setminus K) \right)\cup\left(E_{\alpha}(\mu)\cap
K\right).\]\end{proposition}

\begin{remark}\label{rema}\rm \begin{enumerate} \item If $B=\{0,\cdots,\ell-1\}$,
  $K=[0,1]$ and proposition \ref{ealpha}  is immediate.\item If
  $B=\emptyset$, $K=\emptyset$ and  by proposition \ref{ealpha},  
  $E_{\alpha}(\mu)=E_{\alpha}(\nu)$. In fact, in this case, it is easy to
  prove that    the measures $\mu$ and $\nu$ are strongly equivalent. \item  If $B$ is reduced to a single element, $K$  is a singleton.\item   In all other cases, $K$ is a Cantor set.\end{enumerate}\end{remark}
  
\noindent{\bf Proof of  proposition~\ref{ealpha}.} According to the above  remark, we can suppose that $B\not=\{0,\cdots,\ell-1\}$ and $B\not=\emptyset$. Fix  $x\not\in K$ and
$\alpha>0$. To prove our claim, it is sufficient to show  that $x\in E_\alpha(\mu)$ if and only if
$x\in E_\alpha(\nu)$.
Since $x\not\in K$,  there exists $n(x)\in \N$ and $\epsilon\not\in B$  such that \[\forall  n\ge n(x),\quad I_n(x)=I_{\epsilon_1\cdots\epsilon_{n(x)}\epsilon\cdots\epsilon_n}\] where $(\epsilon_1,\cdots,\epsilon_{n(x)})\in B^{n(x)}$.  With
 obvious  notation, it results from  \eref{mat}
 and \eref{norme} that  \begin{eqnarray*} 
   \fl{{\mu(I_n(x))}\over{\nu(I_n(x))}}={{\Vert M_{I_{n(x)}}M_\epsilon
       M_{x,n}\Vert_1}\over {\Vert M_{I_{n(x)}}M_\epsilon M_{x,n}\Vert}}= {{\left\Vert
         \left(\begin{array}{cc} a(x) & 0\\ b(x)& c(x)\end{array}\right)\left(\begin{array}{cc}
           p_\epsilon &p_{\epsilon+\ell}\\p_{2\ell-1-\epsilon} &
           p_{\ell-1-\epsilon}\end{array}\right) \left(\begin{array}{cc} \tilde a_{x,n} & \tilde b_{x,n}\\ \tilde c_{x,n} &
           \tilde d_{x,n} \end{array}\right)\right\Vert_1}\over {\left\Vert \left(\begin{array}{cc}
           a(x) & 0\\ b(x)& c(x)\end{array}\right)\left(\begin{array}{cc}  p_\epsilon
           &p_{\epsilon+\ell}\\p_{2\ell-1-\epsilon} & p_{\ell-1-\epsilon}
         \end{array}\right) \left(\begin{array}{cc} \tilde a_{x,n} & \tilde b_{x,n}\\ \tilde c_{x,n} & \tilde d_{x,n}
         \end{array}\right)\right\Vert}}\\\ge C(x){{\left\Vert \left(\begin{array}{cc}
           p_\epsilon &p_{\epsilon+\ell}\\p_\epsilon+p_{2\ell-1-\epsilon} &
           p_{\epsilon+\ell}+ p_{\ell-1-\epsilon}\end{array}\right)\left(\begin{array}{cc} \tilde a_{x,n} & \tilde b_{x,n}\\ \tilde c_{x,n} & \tilde d_{x,n} \end{array}\right)\right\Vert_1} \over {\left\Vert\left(\begin{array}{cc}p_\epsilon &p_{\epsilon+\ell}\\p_{\epsilon}+p_{2\ell-1-\epsilon} &p_{\epsilon+\ell}+ p_{\ell-1-\epsilon} \end{array}\right)\left(\begin{array}{cc} \tilde a_{x,n} & \tilde b_{x,n}\\ \tilde c_{x,n} & \tilde d_{x,n} \end{array}\right)\right\Vert}}\\= C(x){{p_\epsilon(\tilde a_{x,n}+\tilde b_{x,n})+p_{\epsilon+\ell}(\tilde c_{x,n}+\tilde d_{x,n})}\over {(2p_\epsilon+p_{2\ell-1-\epsilon})(\tilde a_{x,n}+\tilde b_{x,n})+(2p_{\epsilon+\ell}+p_{\ell-1-\epsilon})(\tilde c_{x,n}+\tilde d_{x,n})}}.\end{eqnarray*}if $n$ is large enough.
Using  $p_\epsilon p_{\epsilon+\ell}>0$, we can find  another constant   $C'(x)>0$ which does not depend on
$n$ such that $\mu(I_n(x))\ge C'(x)\nu(I_n(x)).$    Since 
$\mu\le 2\nu$,  we then  deduce that  \[\left(E_{\alpha}(\mu)\cap
([0,1]\setminus K) \right)=\left(E_{\alpha}(\nu)\cap
([0,1]\setminus K) \right)\] and 
proposition~\ref{ealpha} follows easily.$\square$
\begin{theorem}\label{deco} The  notation are the same as in theorem \ref{tau}.   Suppose  that $B\cap B^*=\emptyset$ where
  $B^*=\{\ell-1-\epsilon,\epsilon\in B\}$.     Then,\[\fl (i)\,\,D_{\mu}=
  D_{\nu}\cup\left[-\log_\ell\left(\max_{i\in B}p_i\right),
    -\log_\ell\left(\min_{i\in B}p_i\right)\right]\quad \mbox{where}\quad D_{\mu}=\{\alpha,\,
E_\alpha(\mu)\not=\emptyset\}.\] (By convention, $\left[-\log_\ell\left(\max_{i\in B}p_i\right),
    -\log_\ell\left(\min_{i\in B}p_i\right)\right]=\emptyset$ if $B=\emptyset$.)   \[\fl(ii)\,\, \forall \alpha\in D_{\mu},\quad \dim(E_\alpha(\mu))=\max\left(\dim(E_\alpha(\nu)),\tilde\tau^*(\alpha)\right).\] \end{theorem}
\noindent{\bf Proof.}  If $B=\emptyset$, the measures $\mu$ and $\nu$ are strongly equivalent and theorem \ref{deco} is easy. 
Now assume that $B\not=\emptyset$. Let  $\pi$ be  the self-similar  measure supported on 
$K$   verifying    \begin{eqnarray}\label{pi} \pi={1\over {\sum_{i\in B} p_i}}\sum_{i\in B} p_i\,\pi\circ S_i^{-1}.\end{eqnarray}
The family $(S_i)_{i\in B}$ satisfying  the Open Set Condition  \cite{hu}, the calculation of the 
   $L^q$-spectrum $\tau_\pi(q)$ is then straightforward~:   
\begin{eqnarray}\label{topi} \forall q\in\R,\quad \tau_\pi(q)=\log_\ell\left(\sum_{i\in B} \left({{p_i}\over {\sum_{i\in B} p_i}}\right)^q\right),\end{eqnarray}(e.g. \cite{cm,o}). 
Moreover,  we have   \begin{eqnarray}\label{ee}\fl\cases{\dim(E_{\alpha}(\pi))=\tau_\pi{}^*(\alpha)&$\mbox{if}\,\,\alpha\in \left[-\log_\ell\left({{\max_{i\in B}(p_i)}\over{\sum_{i\in B} p_i}}\right),-\log_\ell\left({{\min_{i\in B}(p_i)}\over{\sum_{i\in B} p_i}}\right)\right],$\\E_{\alpha}(\pi)=\emptyset & $\mbox{otherwise}.$\\}\end{eqnarray}   
Furthermore, for every    $(\epsilon_1\cdots\epsilon_n)\in B^n$ and   $I=I_{\epsilon_1\cdots\epsilon_n}$, we find that
\[\pi(I)=\left(\sum_{i\in B} p_i\right)^{\log_\ell(\vert I\vert)}\mu(I).\]
 We thus  deduce that 
\[\forall \alpha, \quad  E_{\alpha}(\pi)=E_{\alpha-\log_\ell\left(\sum_{i\in
      B}p_i\right)}(\mu)\cap K,\] or equivalently  \[\forall \alpha,
\quad E_{\alpha+\log_\ell\left(\sum_{i\in B}p_i\right)}(\pi)= E_\alpha(\mu)\cap
K.\]   It follows  from  \eref{ee} that    \begin{eqnarray*} E_\alpha(\mu)\cap
  K\not=\emptyset \Leftrightarrow \alpha\in \left[-\log_\ell\left(\max_{i\in
        B}p_i\right), -\log_\ell\left(\min_{i\in
        B}p_i\right)\right],\end{eqnarray*} 
and, if  $-\log_\ell\left(\max_{i\in B}p_i\right)\le\alpha\le
  -\log_\ell\left(\min_{i\in B}p_i\right)$, then          \begin{eqnarray}\label{deco1}\dim(E_{\alpha}(\mu)\cap K)=\tau_\pi{}^*\left(\alpha+\log_\ell\left(\sum_{i\in B}p_i\right)\right)=\tilde\tau^*(\alpha).\end{eqnarray} 

Proposition \ref{ealpha}  also leads to estimate  $\dim (E_{\alpha}(\nu)\cap
([0,1]\setminus K))$.   The hypothesis   $B\cap B^*=\emptyset$ implies that $T(K)\subset [0,1]\setminus K$. Therefore,
\begin{eqnarray*}\dim(E_{\alpha}(\nu)\cap K)&=&\dim(T(E_{\alpha}(\nu)\cap
K))\\&=&\dim(E_{\alpha}(\nu)\cap T(K))\le \dim(E_{\alpha}(\nu)\cap
([0,1]\setminus K)),\end{eqnarray*}  and we conclude that   \[\dim(E_{\alpha}(\nu))= \dim(E_{\alpha}(\nu)\cap
([0,1]\setminus K)).\]   Theorem \ref{deco} then   follows  from       proposition  \ref{ealpha}  and \eref{deco1}.$\square$  
\begin{remark}\rm \begin{enumerate} \item[1.]  Using  the same ideas, we can also obtain that  \[\forall \alpha,\quad \dim(V_\alpha(\mu))=\max\left(\dim(V_\alpha(\nu)),
 \tilde\tau^*(\alpha)\right),\] where $V_\alpha$ is defined as $E_\alpha$  replacing  $\lim$ by $\liminf$. In other terms,   \[V_\alpha(m)=\,
 \left\{x \in \lbrack 0,1\rbrack,\, \liminf_{n\rightarrow +\infty} -{{\log m
       (I_n(x))}\over {n \log \ell}}=\alpha\right\}.\] 

\item[2.]  Similar results
 can also be established   replacing the Hausdorff dimension $\dim(E_\alpha)$   by the Packing dimension $\Dim(E_\alpha)$. \ \end{enumerate}\end{remark}   
 
We deduce  the following.
\begin{corollary}\label{cor2}  Suppose that  $B\cap B^*=\emptyset$ and that
  $\nu$ satisfies the quasi-Bernoulli property. Then,   \begin{eqnarray*}\fl \forall \alpha,\,\,\dim(E_\alpha(\mu))=\Dim(E_\alpha(\mu))=\dim(V_\alpha(\mu))=\Dim(V_\alpha(\mu))=\max\left(\tau_\nu{}^*(\alpha),
 \tilde\tau^*(\alpha)\right).\end{eqnarray*}\end{corollary}

According to theorem \ref{tau}, the function $\tau_\mu{}^*$ is the Legendre
transform of the maximum of $\tau_\nu$ and $\tilde\tau$. On the other hand, 
by  corollary \ref{cor2},   the dimension of the level sets
$E_\alpha(\mu)$ is given by the maximum of the Legendre transform of
$\tau_\nu$ and   the Legendre transform of $\tilde\tau_\nu$. Since we cannot invert Legendre transform and maximum, we have the following.
\begin{theorem}\label{raj} Suppose that   $B\cap B^*=\emptyset$ and  that
  $\nu$ satisfies the quasi-Bernoulli property. Then, we have the following. 
\begin{enumerate} \item If $\tau_{\mu}{}'(q)$ exists and if 
  $\alpha=-\tau_{\mu}{}'(q)$, then 
  \[\dim(E_{\alpha}(\mu))=\Dim(E_{\alpha}(\mu))=\dim(V_{\alpha}(\mu))=\Dim(V_{\alpha}(\mu))=\tau_\mu{}^*(\alpha).\]\item If
  $\tau_{\mu}{}'(q)$ does not exist and if 
  $-(\tau_{\mu})_+'(q)<\alpha<-(\tau_{\mu})_-'(q),$ then
  \[\dim(E_{\alpha}(\mu))=\Dim(E_{\alpha}(\mu))=\dim(V_{\alpha}(\mu))=\Dim(V_{\alpha}(\mu))<\tau_\mu{}^*(\alpha).\]\end{enumerate} 
  
  Hence, each phase
transition $q$  gives rise  to an  interval $\left(-(\tau_{\mu})_+'(q),-(\tau_{\mu})_-'(q)\right)$ in
 which the multifractal formalism breaks down.\end{theorem}
\begin{remarke}\rm  It is possible to prove that  $\nu$ satisfies the
  quasi-Bernoulli property if and only~if  \begin{eqnarray}\label{nuqb}\mbox{either}\quad\forall\, i\in B,\,\,\,p_i<p_{\ell-1-i}, \quad \mbox{or}\quad   \forall i\in B,\,\,\, p_i>p_{\ell-1-i}.\end{eqnarray}  The multifractal formalism  fails for the measure $\mu$ only in  the  first case.
 Indeed, if for every $ i\in B,\,\,\, p_i>p_{\ell-1-i}$, it is easy to check
that  the measures $\mu$ and $\nu$ are strongly equivalent.\end{remarke} 
\section{Examples}
In this section we construct measures with non-differentiable $L^q$-spectra $\tau(q)$ for which previous results apply. Furthermore, based on these examples, we point out  
  new phenomena in  the multifractal structure
 of self-similar measures.  
\subsection{An  isolated point in the set of  local dimensions}
\begin{figure}[h!]
\begin{center}\subfigure[]{\epsfig{figure=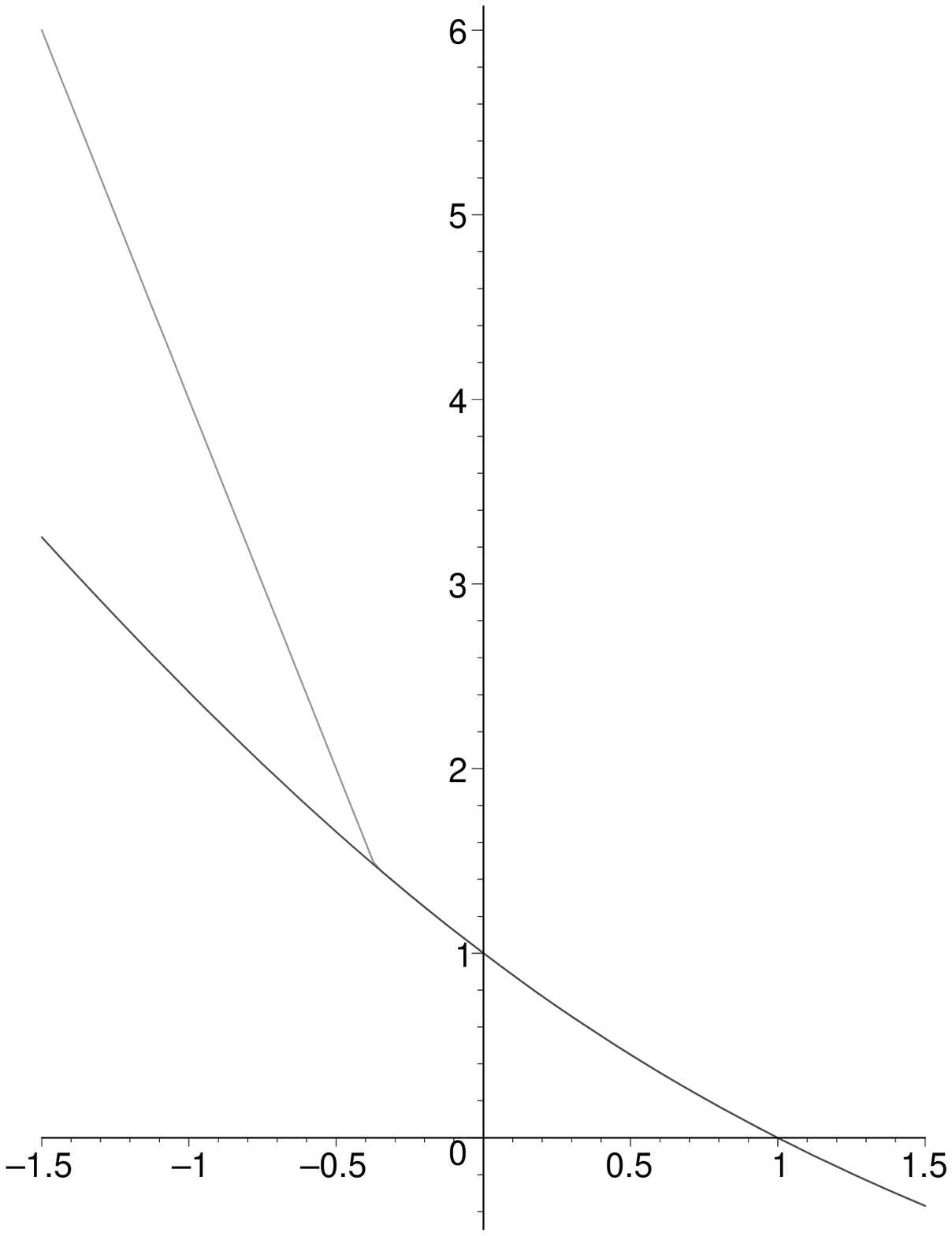,width=4.5cm,height=4cm,angle=360}}\quad\quad
  \subfigure[]{\epsfig{figure=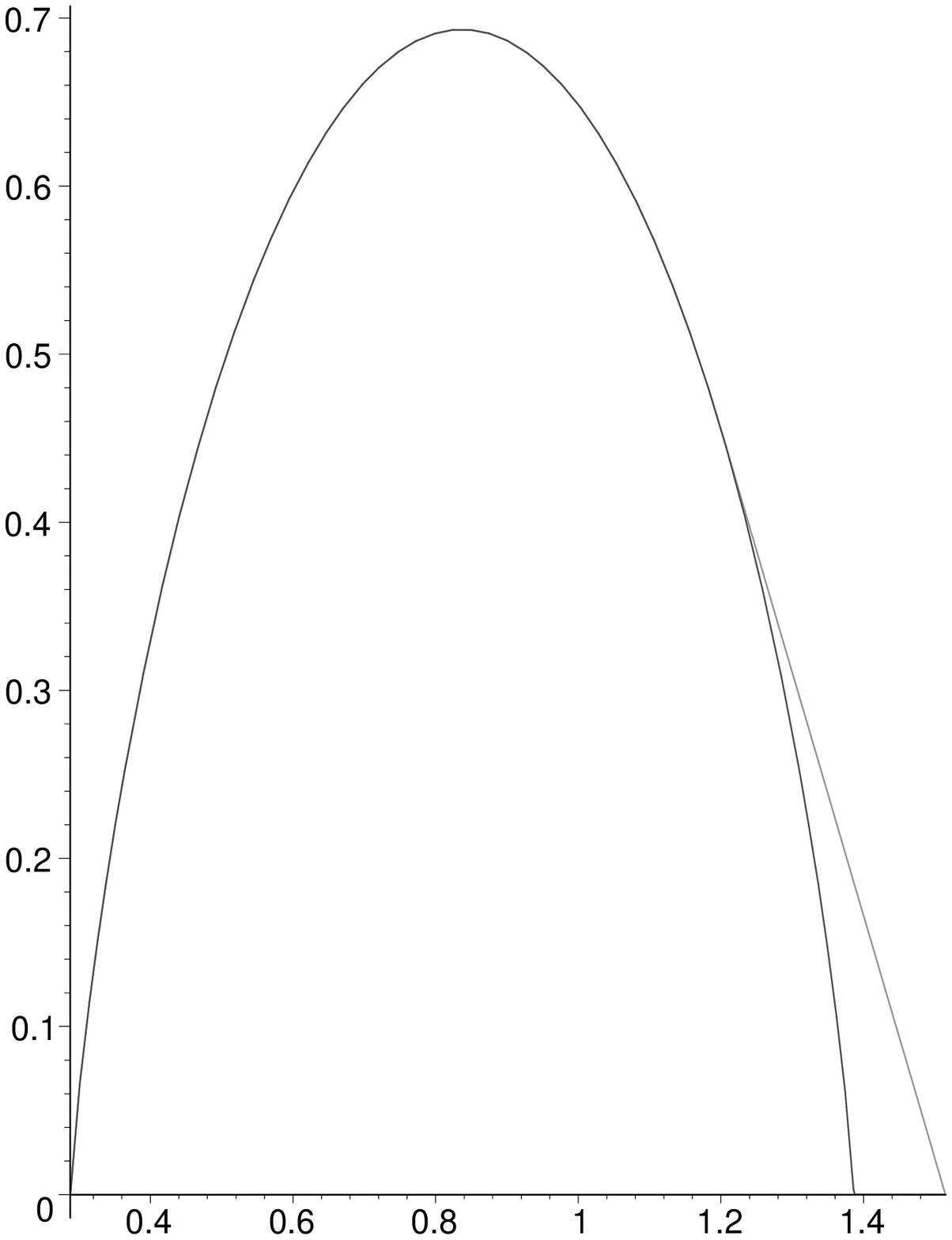,width=4.5cm,height=4cm,angle=360}}
\caption{ (a) $\tau_{\mu}$ is not differentiable. (b) The singularity spectrum of $\mu$, given by $\tau_{\nu}{}^*$, differs from $\tau_\mu{}^*$.}\label{figu2}\end{center}
\end{figure}

Let us  take $\ell=2$ and  consider the  probability measure $\mu$ verifying 
\begin{eqnarray}\label{l=2}\mu=p_0\mu\circ S_0^{-1}+p_1\mu\circ S_1^{-1}+p_2\mu\circ S_2^{-1},\end{eqnarray}
where $S_0(x)=x/2$, $S_1(x)=x/2+1/2$ and  $S_2(x)=-x/2+1/2$. We assume that  $p_0 p_1 p_2
>0$ and  $p_1<p_0$. With the notation previously introduced, we have $B=\{1\}$  and  $K=\{1\}$. Moreover, by theorem~\ref{tau}, \[\forall q\in\R,\quad\tau_{\mu}(q)=\max(\tau_{\nu}(q),q\log_2(p_1)).\]
Thus, in order to get a phase transition for the function $\tau_{\mu}$, we
have to compare $\tau_{\nu}{}'(-\infty)$ and  $\log_2(p_1)$. For every $I\in \mathcal{F}_n$, by iterating 
 \eref{l=2}, we get $\nu(I)\ge  (p_{-})^n$ where    $p_-  =\min(p_0,p_1+p_2)$. We   easily  deduce that 
$-\tau_{\nu}{}'(-\infty)\le-\log_2( p_-)
  <-\log_2(p_1)$. Since   $\tau_{\nu}(q)\ge q\log_2(p_1)$ for $q=0$, we conclude that there exists  $q_0<0$ such that 
 \begin{equation}\tau_{\mu}(q)=\cases{ q\log_2(p_1) & $\mbox{ if }\quad  q\le q_0$,\\ \tau_{\nu}(q) & $\mbox{ if }\quad q \ge q_0$,\\}\end{equation}
  and the $L^q$-spectrum $\tau_{\mu}(q)$ is not differentiable at  $q=q_0$
(see  \fref{figu2}(a)). 
Furthermore, by \eref{nuqb}, $\nu$ satisfies the quasi-Bernoulli
property. Using theorem \ref{deco} and corollary
\ref{cor2} we  deduce the following    .  
 \begin{theorem}\label{2bis} Let $\mu$ be the measure satisfying
  \eref{l=2}. Then, \[D_{\mu}=D_{\nu}\cup \{-\log_2(p_1)\}=
\left(-\tau_{\nu}{}'(+\infty),-\tau_{\nu}{}'(-\infty)\right)  \cup
\{-\log_2(p_1)\},\] and \begin{eqnarray*} \dim(E_\alpha(\mu))=\cases{
     \tau_{\nu}{}^*(\alpha) & $\mbox{if } \,\,
     \alpha \in
    \left(-\tau_{\nu}{}'(+\infty),-\tau_{\nu}{}'(-\infty)\right),$\\0 &$\mbox{if }\,\,\alpha=-\log_2(p_1).$\\}\end{eqnarray*}\end{theorem}
\begin{remark}\rm \begin{enumerate} \item[1.]  Since $-\tau_{\nu}{}'(-\infty) <-\log_2(p_1)$,  $D_\mu$   contains an isolated point.  In this sense, the situation is close to  the ones  obtained for the Erd\"os measure  
  and  for  the $3$-time convolution of the Cantor measure (e.g \cite{fl1, hl}).  Note that in our situation,  the  value of  $[p_0, p_1, p_2]$ is not a matter.

\item[2.]  It is easy to show that  \begin{eqnarray*}\fl
    \tau_\mu{}^*(\alpha)=\cases{\tau_{\nu}{}^*(\alpha)  & $\mbox{ if
      }\quad -\tau_{\nu}{}'(+\infty)\le\alpha\le -\tau_{\nu}{}'(q_0),$\\
      {{\tau_\nu{}^*(q_0)}\over {\log_2(p_1)-\tau_\nu{}'(q_0)}}\left(\alpha+\log_2(p_1)\right)  & $\mbox{ if }\quad    -\tau_{\nu}{}'(q_0)\le \alpha \le
      -\log_2(p_1).$\\}\end{eqnarray*} Thus,  
\begin{eqnarray*}\fl\cases{\forall\alpha  \in   (-\tau_{\nu}{}'(+\infty)
    ,-\tau_{\nu}{}'(q_0)\rbrack, \quad \dim(E_{\alpha}(\mu))=\dim(V_{\alpha}(\mu))=\tau_{\nu}{}^*(\alpha)=\tau_\mu{}^*(\alpha),&\\
\forall \alpha \in   ( -\tau_{\nu}{}'(q_0),-\tau_{\nu}{}'(-\infty)
\rbrack,\quad\dim(E_{\alpha}(\mu))=\dim(V_{\alpha}(\mu))=\tau_{\nu}{}^*(\alpha)<\tau_\mu{}^*(\alpha).&\\}\end{eqnarray*}
 Contrary to the usual  situation,  the singularity  spectrum of
  $\mu$ is not     given by the Legendre
transform of  $\tau_{\mu}$  but instead  by the Legendre
transform of  an auxiliary function.     \Fref{figu2}(b) illustrates  this  phenomenon.
\item[3.]  The measure $\mu$ may be used   to estimate the Hausdorff dimension of 
  self-affine graphs studied by  McMullen~\cite{mc}, Prsytycki and Urba\'nski
  \cite{pu,ur}. More details can be found in    \cite{t1,t2}.\end{enumerate} 
\end{remark}
\subsection{Non-concave spectra}
Subsection 3.1 and several papers   deal with    measures for which the 
$L^q$-spectrum  $\tau(q)$ is     not differentiable at a single point  $q=q_0$ and is linear
for  $q\le
q_0$ (e.g.  \cite{f,fl1,hl,ln,ost}). In this part we construct measures with non-differentiable and
strictly concave
$L^q$-spectra.   That leads to new  situations for the
multifractal analysis of self-similar  measures.

Let us   take $\ell=4$ and   consider  the probability measure $\mu$ 
satisfying 
 \begin{eqnarray}\label{kp2} \mu =\sum_{i=0}^5 p_i\, \mu \circ S_i^{-1},\end{eqnarray}where
\[S_0(x)={x\over 4} ,\quad S_1(x)={x\over 4}+{1\over 4},\quad S_2(x)={ x\over
  4}+{1\over 2},\] \[S_3(x)={x\over 4}+{3\over 4},   \quad S_4(x)=-{x\over
  4}+{1\over 4} \quad\mbox{and}\quad S_5 (x)=-{x\over 4}+{1\over 2}.\]
In this case,  $B=\{2,3\}$ and  $K$ is the Cantor set  whose  points only contains digits $2$ and
$3$ in  their  $4$-adic expression,   i.e.  $K=\left\{x=\sum \epsilon_i/4^i, \, \epsilon_i=2 \,\, \mbox{or} \,\, 3,\, \forall i\in\N^*\right\}$. 

\begin{figure}[h!]
\begin{center}
\subfigure[]{\epsfig{figure=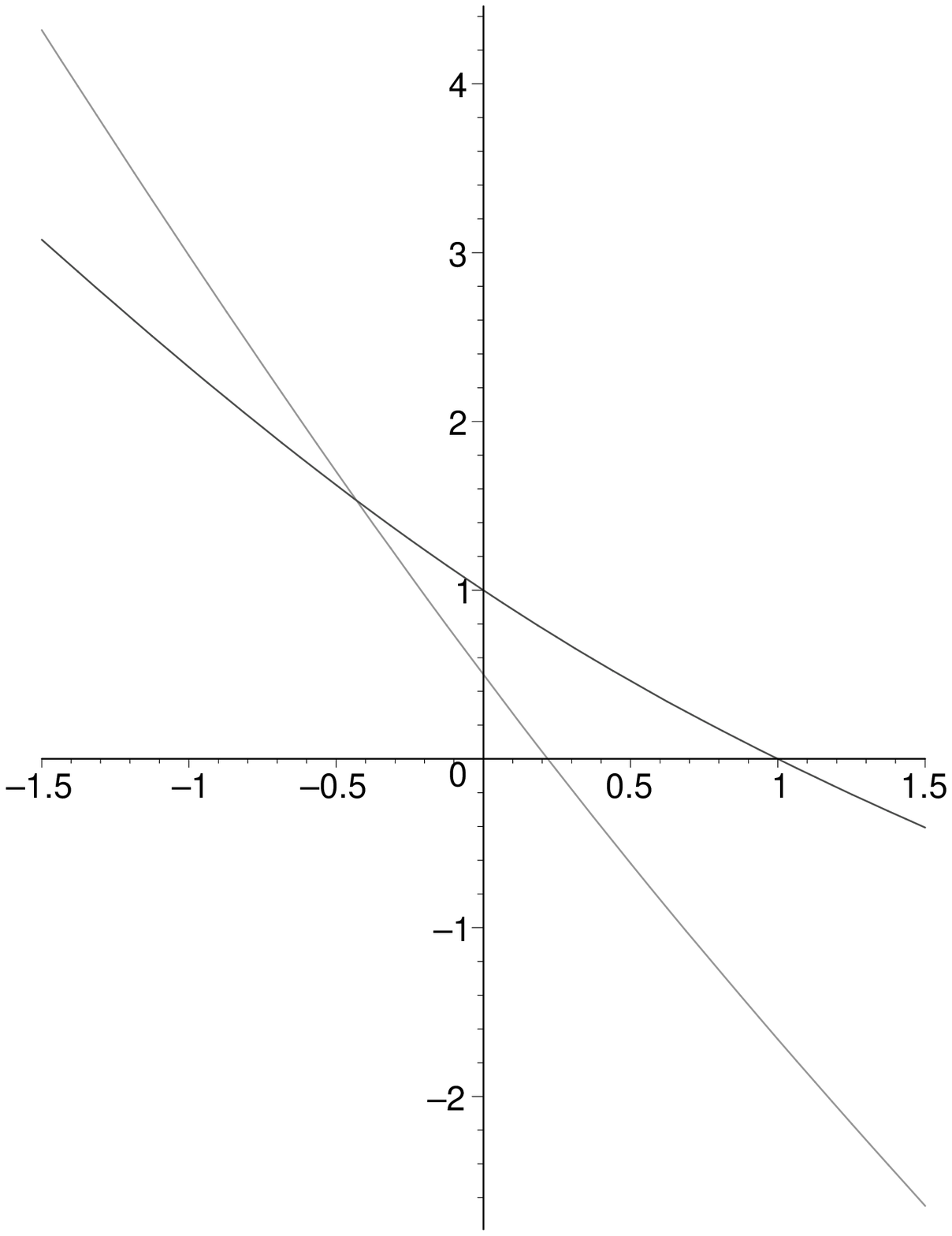,width=4.5cm,height=4cm,angle=360}}\quad\quad
  \subfigure[]{\epsfig{figure=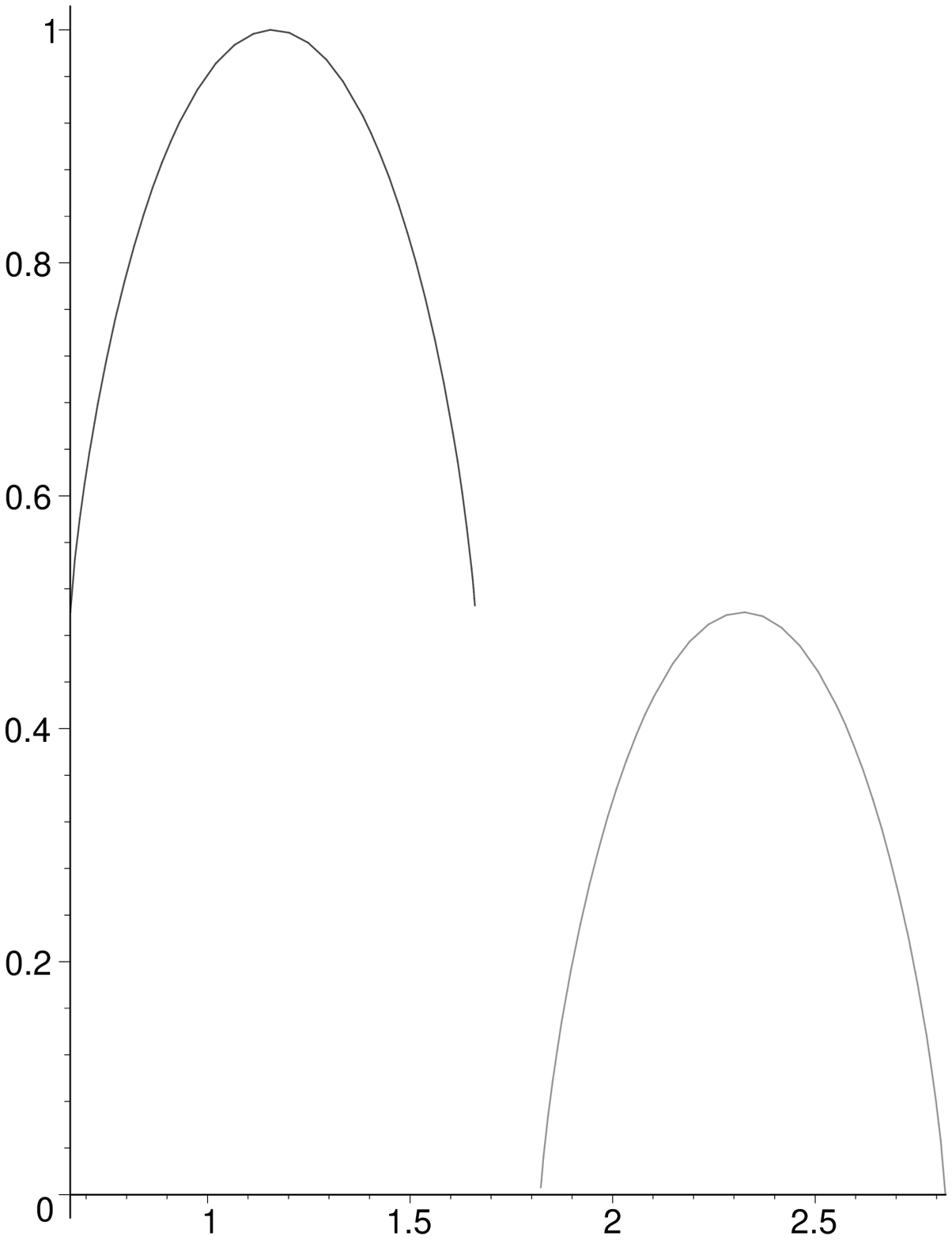,width=5.0cm,height=4cm,angle=360}}
\caption{ (a) $\tau_{\mu}$ is not differentiable. (b) The singularity spectrum of
  $\mu$ is supported by a union of two disjoint
  intervals.}\label{fig5}\end{center}\end{figure}

By theorem~\ref{tau}, $\tau_\mu(q)=\max(\tau_\nu(q),\log_4(p_2{}^q+p_3{}^q))$. In order to compute  $\tau_{\nu}$, we assume that the $p_i{}'s$ verify $p_0=p_3+p_4$ and $ p_1=p_2+p_5.$  In this situation, it is easy to show that $\nu$ is a  multinomial  measure (see \cite{t1}). The calculation of  $\tau_\nu$ is then   straightforward~:    $\tau_{\nu}(q)=1/2
    +\log_4(p_0{}^q+p_1{}^q)$. Therefore, there exists $q_0<0$ such that
 \begin{eqnarray}\label{res} \tau_{\mu}(q)=\cases{\log_4(p_2{}^q+p_3{}^q) &$\mbox{ if }\quad  q\le q_0$,\\ {1\over 2}+\log_4(p_0{}^q+p_1{}^q) & $\mbox{ if
     }\quad q \ge q_0$,\\}\end{eqnarray}   and  
 $\tau_{\mu}(q)$ is not differentiable at   $q=q_0$ (see \fref{fig5}(a)).  

 Moreover, if we  denote   $p_0\vee p_1$ ($p_0\wedge p_1$)    the maximum (minimum)  of
$p_0$ and $p_1$,  we get  $D_{\nu}=[-\log_4(p_0\vee p_1),-\log_4(p_0\wedge p_1)]$ and 
  $\dim(E_\alpha(\nu))=\tau_{\nu}{}^*(\alpha)\ge 1/2$, for all~$\alpha\in D_{\nu}$.    It follows  from theorem \ref{deco} that   \[D_{\mu}=[-\log_4(p_0\vee p_1),-\log_4(p_0\wedge p_1)]\cup[-\log_4(p_2\vee p_3),-\log_4(p_2\wedge p_3)],\] and   \begin{eqnarray*} \dim(E_\alpha(\mu))=\cases{
     \tau_{\nu}{}^*(\alpha) & $\mbox{if } \,\,
     -\log_4(p_0)\le\alpha\le -\log_4(p_1),$\\\tilde\tau^*(\alpha) &$\mbox{if }\,\, -\log_4(p_2\vee p_3)\le\alpha\le-\log_4(p_2\wedge p_3).$\\}\end{eqnarray*}
Thus,  if $p_3<p_1\le p_0$,  the singularity spectrum of $\mu$ is supported by a  union of mutually  disjoint intervals  and differs from $\tau_\mu^*(\alpha)$  for  $-(\tau_{\mu})_+'(q_0)<\alpha<-(\tau_{\mu})_-'(q_0)$ (see \fref{fig5}(b)). To the  best of our  knowledge,   self-similar measures with such multifractal structures have not previously   appeared  in the litterature.

\subsection{Two  phase transitions}

Until now  we have studied measures for which  the $L^q$-spectrum $\tau(q)$ is
not differentiable at  one single   point  $q_0<0$. In this part we propose examples  with two phase transitions. Let  us take $\ell=5$ and  consider the probability measure $\mu$ satisfying  \begin{eqnarray}\label{kp3} \mu =\sum_{i=0}^7 p_i\, \mu \circ S_i^{-1},\end{eqnarray}
where
\[\fl S_0(x)={x\over 5} ,\quad S_1(x)={x\over 5}+{1\over 5},\quad S_2(x)={ x\over
  5}+{2\over 5}, \quad S_3(x)={ x\over 5}+{3\over 5} \] \[\fl S_4(x)={x\over 5}+{4\over 5},   \quad S_5(x)=-{x\over 5}+{1\over 5}, \quad S_6(x)=-{x\over 5}+{2\over 5}\quad\mbox{and}\quad   S_7(x)=-{x\over 5}+{3\over 5}.\]
In this case,   $B=\{3,4\}$ and  $K=\left\{x=\sum \epsilon_i/5^i, \, \epsilon_i=3 \,\, \mbox{or} \,\, 4,\, \forall i\in\N^*\right\}$.    We
suppose that the coefficients~$p_i{}'s$  verify   $p_0=p_4+p_5$, $p_1=p_3+ p_6$ and   $p_2=p_7$. As in section 6.2, we get       $\tau_{\mu}(q)=\max(\log_5(2p_0{}^q+2p_1{}^q+(2p_2)^q),\log_5(p_3{}^q+p_4{}^q)),$ and \[\fl\forall \alpha\in
  D_{\nu}=[-\log_5(p_0\vee p_1\vee 2p_2),-\log_4(p_0\wedge p_1\wedge 2p_2)],
  \,\,\, \dim(E_\alpha(\nu))=\tau_{\nu}{}^*(\alpha).\]
  In order to  have  $\tau_{\mu}(q)=\tau_{\nu}(q)$ for  large negative
 $q$, we choose  $p_2$ sufficiently small. For  example, if we take  
$p_0=0.35$,  $ p_1=0.14,$   $p_2=0.01,$  $p_3=0.03$ 
and $p_4=0.025$, the  equation $\tau_{\nu}( q)=\tilde\tau(q)$
has  two solutions $q_0$ and $q_1$ corresponding to the points of non-differentiability of  $\tau_{\mu}(q)$. By  theorem \ref{deco},    $D_{\mu}=[-\log_5(p_0),-\log_5(2p_2)]$ and \ \begin{eqnarray*} \dim(E_\alpha(\mu))=\cases{\tau_{\nu}{}^*(\alpha) &  if $\quad -\log_5(p_0)\le\alpha\le \alpha_0$,\\
\tilde\tau^*(\alpha) & if $\quad \alpha_0\le\alpha\le\alpha_1$,\\
\tau_{\nu}{}^*(\alpha) & if $\quad\alpha_1\le \alpha\le
-\log_5(p_2),$\\}\end{eqnarray*} where     $\alpha_0$ and $\alpha_1$ denote  the solutions of the equation $\tau_\nu{}^*(\alpha)=\tilde\tau^*(\alpha)$. From the expression of the Legendre transform of 
$\tau_{\mu}=\max(\tau_{\nu},\tilde\tau)$, it follows that     \[\fl\forall\alpha\in\left(-(\tau_{\mu})_+'(q_0),-(\tau_{\mu})_-'(q_0)\right)\cup \left(-(\tau_{\mu})_+'(q_1),-(\tau_{\mu})_-'(q_1)\right),\quad
\dim(E_{\alpha}(\mu))<\tau_\mu{}^*(\alpha).\]
\begin{figure}[h!]
\begin{center}\subfigure[]{\epsfig{figure=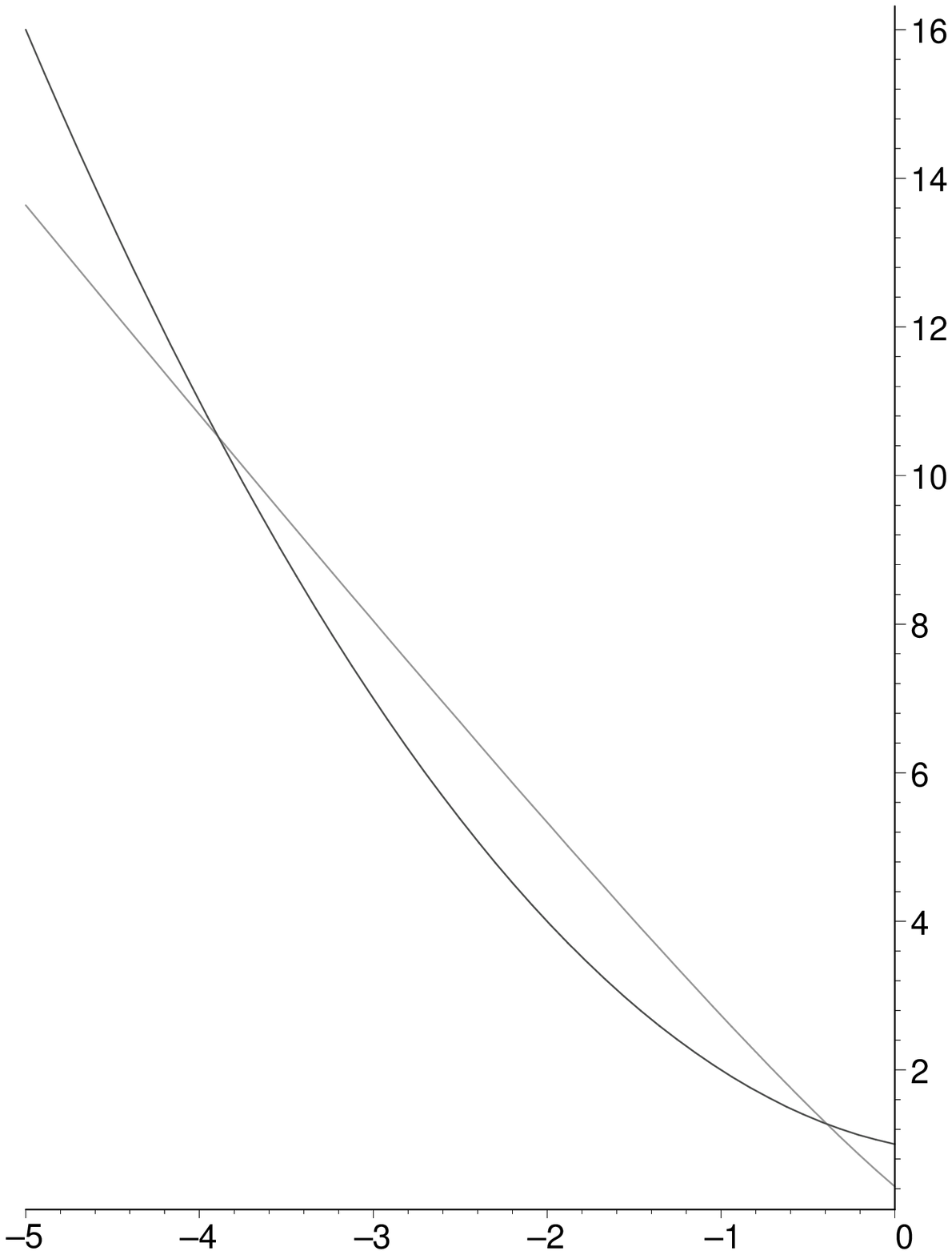,width=4.5cm,height=3.7cm,angle=360}}\quad\quad
  \subfigure[]{\epsfig{figure=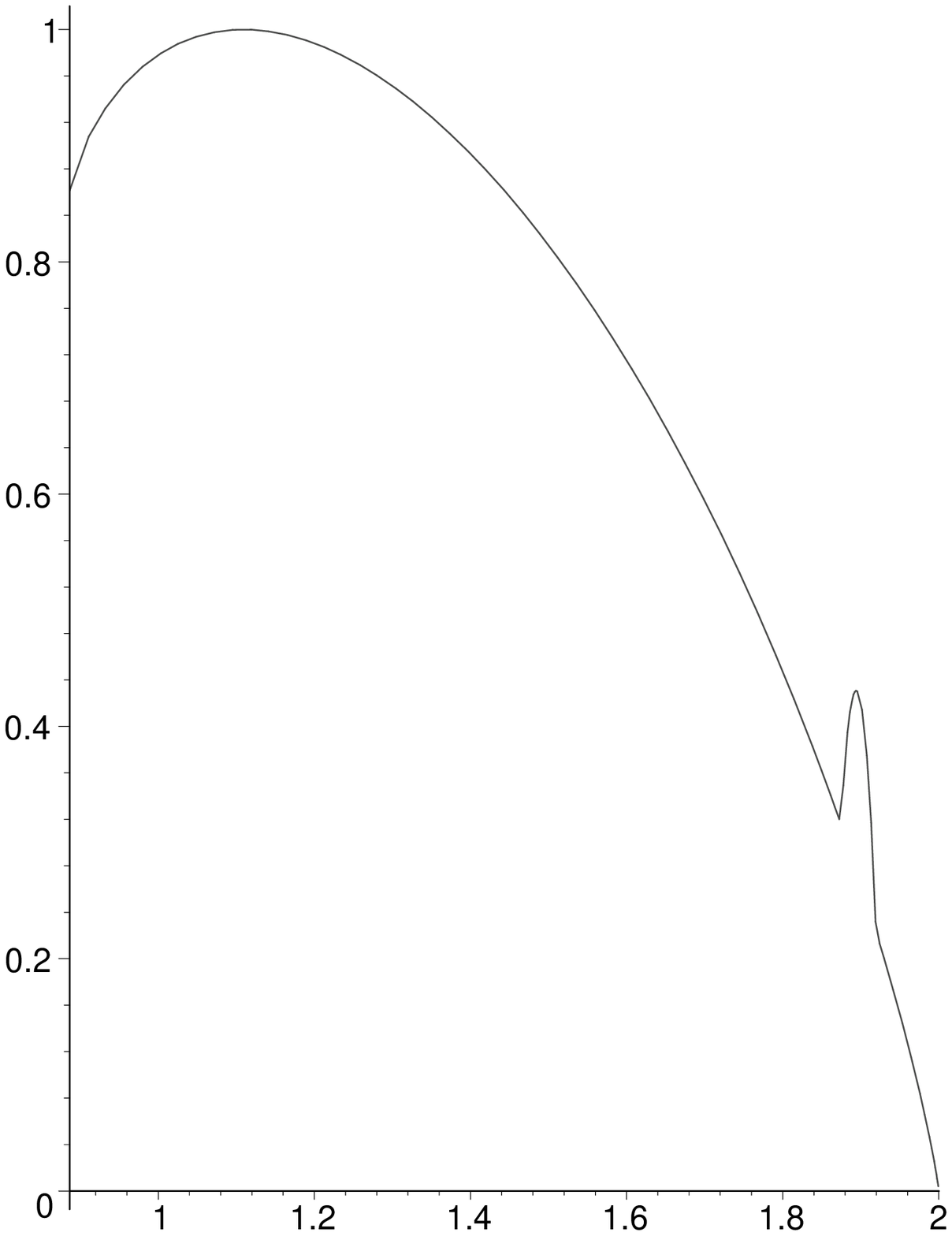,width=5.0cm,height=3.7cm,angle=360}}
\caption{$\tau_{\mu}=\max(\tau_{\nu}, \tilde\tau)$ is not differentiable at two points. The
  singularity  spectrum of $\mu$ is not  concave  and differs from
  $\tau_\mu{}^*$.}\label{fig14} 
\end{center}
\end{figure}
%\newpage  
\subsection{More phase transitions}
In this part  we describe a way  to construct measures with
an arbitrarily large number $N$  of phase transitions. Theorem \ref{tau} leads us to  find conditions on the 
 $p_i{}'s$ such  that the equation $\tau_\nu(q)=\tilde\tau(q)$ has $N$
solutions. Since $\tau_\nu(0)\ge \tilde\tau(0)$, we have to distinguish the case   where $N$ is  odd from the  case where   $N$  is even. 

First, assume that $N$ is odd.  Let us  take $\ell=2N$,  
$B=\{\ell/2,\cdots,\ell-1\}$ and 
 suppose  that $
p_i=p_{i+\ell}+p_{\ell-1-i}$, for all $0\le i\le\ell/2-1$. In this case, the arguments developed  in section 6.2 imply that 
\[\tau_\mu(q)=\max(\tau_\nu(q),\tilde\tau(q))=\max\left(\log_{2N}\left(\sum_{i=0}^{N-1}2p_i{}^q\right),\log_{2N}\left(\sum_{i=N}^{2N-1}p_i{}^q\right)\right).\] Moreover, since $\ell=2N$, we can choose the   $p_i{}'s$ such that the equation
\[\sum_{i=0}^{N-1}2p_i{}^q  =  \sum_{i=N}^{2N-1}p_i{}^q\] has $N$ solutions. These solutions correspond to  the    phase transitions for     the $L^q$-spectrum $\tau_\mu(q)$.% (see \fref{arb}(a) in the case $N=5$). 

Assume now that $N$ is even. To ensure  that  $\tau_\mu(q)=\tau_\nu(q)$ for large negative~$q$, tools  used in section 6.3 suggest to take $\ell$ odd. Let  $\ell=2N+1$ and   $B=\{N+1,\cdots,2N\}$. Under the conditions,  for all $0\le i\le N-1$,   $
p_i=p_{i+\ell}+p_{\ell-1-i}$  and   $p_N=p_{N+\ell}$, we get  \[\tau_\mu(q)=\max\left(\log_{2N+1}\left(\sum_{i=0}^{N-1}2p_i{}^q+(2p_N)^q\right),\log_{2N+1}\left(\sum_{i=N+1}^{2N}p_i{}^q\right)\right).\] Thus, in order to have $\tau_\mu(q)= \tau_\nu(q)$ for large negative $q$,  we also suppose that  $2p_N<\min(p_i,\, N+1\le i \le 2N).$ Once again, we can choose the  $p_i{}'s$  such that the equation
\[\sum_{i=0}^{N-1}2p_i{}^q +(2p_N)^q =  \sum_{i=N+1}^{2N}p_i{}^q\] has $N$
solutions. They correspond   to the   
phase transitions for  the function $\tau_\mu$.

 More details about these examples can be found in \cite{t1}. %(see \fref{arb}(b) in the case $N=4$).   

%\begin{figure}[h!]\begin{center}\subfigure[]{\epsfig{figure=figure4a.eps ,width=4.5cm,height=3.3cm,angle=360}}\quad
  %\subfigure[]{\epsfig{figure=figure4b.eps,width=4.5cm,height=3.3cm,angle=360}}\caption{The graph of the function
 % $\tau_\nu-\tilde\tau$. On the left, $\ell=10$ and
  %$p_5=0.00482<p_4=0.5-\sum_{i=0}^3 p_i$ $<p_6=p_7=0.00492<p_3=0.00516$ $<p_8=p_{9}=0.0058<p_0=p_1=p_2=0.49/3.$ On the right, $\ell=9$ and $p_4=0.003<p_5=p_6=0.0062<$$p_3=0.007<p_7=p_8=0.01<p_0=p_1=p_2=0.49/3$. The zeros of   $\tau_\nu-\tilde\tau$  correspond to the phase transitions of $\tau_\mu$.}
%\label{arb} \end{center}
%\end{figure}

\ack 
The author had a position in Clermont-Ferrand (Universit\'e Blaise Pascal) when this work was achevied and he thanks this institution. He also  would like to thank Yanick Heurteaux and Anathassios Batakis    for reading the manuscript
carefully and suggesting some improvements.

\end{document}